\documentstyle{article}
\begin{document}
\setlength{\textheight}{574pt}
\setlength{\textwidth}{432pt}
\setlength{\oddsidemargin}{18pt}
\setlength{\topmargin}{14pt}
\setlength{\evensidemargin}{18pt}
\newtheorem{theorem}{Theorem}[section]
\newtheorem{lemma}{Lemma}[section]
\newtheorem{corollary}{Corollary}[section]
\newtheorem{conjecture}{Conjecture}
\newtheorem{remark}{Remark}[section]
\newtheorem{definition}{Definition}[section]
\newtheorem{problem}{Problem}
\newtheorem{proposition}{Proposition}[section]
\title{{\bf SUBADJUNCTION THEOREM FOR PLURICANONICAL DIVISORS}}
\date{November, 2001}
\author{Hajime TSUJI}
\maketitle
\begin{abstract}
We prove a subadjunction theorem 
which relates the multi-adjoint linear system of 
the ambient space and the linear system of 
the restricted bundle on a subvariety. 
MSC 32H25 
\end{abstract}
\section{Introduction}

Let $M$ be a complex manifold and  $L$ be a line bundle  on $M$
and $S$ be a submanifold of $M$. 
It is a basic question whether the restriction map 
\[
H^{0}(M,{\cal O}_{M}(L)) 
\rightarrow H^{0}(S,{\cal O}_{S}(L))
\]
is surjective. 

In this paper we shall consider this question for 
multi-adjoint type line bundles under certain 
geometric conditions. 

Let us state our result precisely. 
 Let $M$ be a complex manifold of dimension $n$ and let $S$ be a closed complex submanifold of $M$. 
Then we consider a class of continuous function $\Psi : M\longrightarrow [-\infty , 0)$  such that  
\begin{enumerate}
\item $\Psi^{-1}(-\infty ) \supset S$,
\item if $S$ is $k$-dimensional around a point $x$, there exists a local 
coorinate $(z_{1},\ldots ,z_{n})$ on a neighbourhood of $x$ such that 
$z_{k+1} = \cdots = z_{n} = 0$ on $S\cap U$ and 
\[
\sup_{U\backslash S}\mid \Psi (z)-(n-k)\log\sum_{j=k+1}^{n}\mid z_{j}\mid^{2}\mid < \infty .
\]
\end{enumerate} 
The set of such functions $\Psi$ will be denoted by $\sharp (S)$. 

For each $\Psi \in \sharp (S)$, one can associate a positive measure 
$dV_{M}[\Psi ]$ on $S$ as the minimum element of the 
partially ordered set of positive measueres $d\mu$ 
satisfying 
\[
\int_{S_{k}}fd\mu \geq 
\limsup_{t\rightarrow\infty}\frac{2(n-k)}{\sigma_{2n-2k-1}}
\int_{M}f\cdot e^{-\Psi}\cdot \chi_{R(\Psi ,t)}dV_{M}
\]
for any nonnegative continuous function $f$ with 
$\mbox{supp}\, f\subset\subset M$.
Here $S_{k}$ denotes the $k$-dimensional component of $S$,
$\sigma_{m}$ denotes the volume of the unit sphere 
in $\mbox{\bf R}^{m+1}$, and 
$\chi_{R(\Psi ,t)}$ denotes the characteristic funciton of the set 
\[
R(\Psi ,t) = \{ x\in M\mid -t-1 < \Psi (x) < -t\} .
\]

\begin{theorem}
Let $M$ be a complex manifold with a continuous volume form $dV_{M}$,
let $L$ be a holomorphic line bundle over $M$ with a $C^{\infty}$-hermitian metric  $h_{L}$, let $S$ be a compact complex submanifold of $M$,
let $\Psi : M \longrightarrow [-\infty ,0)$ and let $K_{M}$ be the canonical bundle of $M$.
Suppose that the followings are satisfied. 
\begin{enumerate}
\item There exists a closed set $X\subset M$ such that 
\begin{enumerate}
\item $X$ is locally negligble with respect to $L^{2}$-holomorphic functions, i.e., 
for any local coordinate neighbourhood $U\subset M$ and for any $L^{2}$-holomorphic function $f$ on $U\backslash X$, there exists a holomorphic function 
$\tilde{f}$ on $U$ such that $\tilde{f}\mid U\backslash X = f$.
\item $M\backslash X$ is a Stein manifold which intersects with every component of $S$.  
\end{enumerate}
\item $\Psi \in  \sharp (S) \cap C^{\infty}(M\backslash S)$, 
\item $\Theta_{h\cdot e^{-(1+\epsilon )\Psi}}\geq 0$ for 
every $\epsilon \in [0,\delta ]$ for some $\delta > 0$,
\item there is a positive line bundle on $M$.
\end{enumerate}
Then every element of  $H^{0}(S,{\cal O}_{S}(m(K_{M}+L)))$ extends to an element of 
$H^{0}(M,{\cal O}_{M}(m(K_{M}+L)))$. 
\end{theorem}

One may think that the assumption on the existence of
the function $\Psi$ is somewhat technical or restrictive. 
But as one see in the last section, this is not the case. 
In fact one may construct such a function by using 
an effective {\bf Q}-divisor on  $M$.
In the last section we shall formulate a variant of 
Theorem 1.1 which would be more useful but a little bit 
more complicated to formulate. 

This paper is a continuation of \cite{t}
and the most of the ideas are transplanted from \cite{t}. 
The results in this paper may be considered as 
a generalization of \cite{t} to the case of nontrivial normal bundles.
We also note that there exists another type of subadjunction theorem 
due to Y. Kawamata (\cite{ka}). 

The key point of the proof is the extension of closed positive $(1,1)$-current 
by using an algebraic approximation.  
\section{Preliminaries}
\subsection{$L^{2}$-extension theorem} 

Let $M$ be a complex manifold and let $(E,h)$ be a holomorphic hermitian vector 
bundle over $M$. 
Given a positive measure $d\mu_{M}$ on $M$,
we shall denote $A^{2}(M,E,h,d\mu_{M})$ the space of 
$L^{2}$ holomorphic sections of $E$ over $M$ with respect to $h$ and 
$d\mu_{M}$. 
Let $S$ be a closed  complex submanifold of $M$ and let $d\mu_{S}$ 
be a positive measure on $S$. 
The measured submanifold $(S,d\mu_{S})$ is said to be a set of 
interpolation for $(E,h,d\mu_{M})$, or for the 
sapce $A^{2}(M,E,h,d\mu_{M})$, if there exists a bounded linear operator
\[
I : A^{2}(S,E\mid S,h,d\mu_{S})\longrightarrow A^{2}(M,E,h,d\mu_{M})
\]
such that $I(f)\mid S = f$ for any $f$. 
$I$ is called an interpolation operator.
Let $M$ be a complex manifold and let $S$ be a closed 
complex submanifold of $M$.
Let $(L,h)$ be a singular hermitian line bundle on 
$M$. 

Let $dV$ be any continuous nowhere degenerate volume form 
on $M$. 
The following theorem is crucial for the proof of Theorem 1.1.

\begin{theorem}(\cite[Theorem 4]{o})
Let $M$ be a complex manifold with a continuous volume form $dV_{M}$,
let $E$ be a holomorphic vector bundle over $M$ with $C^{\infty}$-fiber 
metric $h$, let $S$ be a closed complex submanifold of $M$,
let $\Psi\in \sharp (S)$ and let $K_{M}$ be the canonical bundle of $M$.
Then $(S,dV_{M}(\Psi ))$ is a set of interpolation 
for $(E\otimes K_{M},h\otimes (dV_{M})^{-1},dV_{M})$, if 
the followings are satisfied.
\begin{enumerate}
\item There exists a closed set $X\subset M$ such that 
\begin{enumerate}
\item $X$ is locally negligble with respect to $L^{2}$-holomorphic functions, i.e., 
for any local coordinate neighbourhood $U\subset M$ and for any $L^{2}$-holomorphic function $f$ on $U\backslash X$, there exists a holomorphic function 
$\tilde{f}$ on $U$ such that $\tilde{f}\mid U\backslash X = f$.
\item $M\backslash X$ is a Stein manifold which intersects with every component of $S$. 
\end{enumerate}
\item $\Theta_{h}\geq 0$ in the sense of Nakano,
\item $\Psi \in \sharp (S)\cap C^{\infty}(M\backslash S)$,
\item $e^{-(1+\epsilon )\Psi}\cdot h$ has semipositive 
curvature in the sense of Nakano for every $\epsilon \in [0,\delta]$ 
for some $\delta > 0$.
\end{enumerate}
Under these conditions, there exits a constant $C$ and an interpolation operator 
from $A^{2}(S,E\otimes K_{M}\mid S,h\otimes (dV_{M})^{-1}\mid S,dV_{M}[\Psi ])$
to $A^{2}(M,E\otimes K_{M},h\otimes (dV_{M})^{-1}.dV_{M})$ whose 
norm does not exceed $C\delta^{-3/2}$.
If $\Psi$ is plurisubharmonic, the interpolation operator can be chosen 
so that its norm is less than $2^{4}\pi^{1/2}$.
\end{theorem}
\subsection{Analytic Zariski decompositon}
In this subsection we shall introduce the notion of analytic Zariski decompositions. 
By using analytic Zariski decompositions, we can handle  big line bundles
like  nef and big line bundles.
\begin{definition}
Let $M$ be a compact complex manifold and let $L$ be a holomorphic line bundle
on $M$.  A singular hermitian metric $h$ on $L$ is said to be 
an analytic Zariski decomposition, if the followings hold.
\begin{enumerate}
\item $\Theta_{h}$ is a closed positive current,
\item for every $m\geq 0$, the natural inclusion
\[
H^{0}(M,{\cal O}_{M}(mL)\otimes{\cal I}(h^{m}))\rightarrow
H^{0}(M,{\cal O}_{M}(mL))
\]
is an isomorphim.
\end{enumerate}
\end{definition}
\begin{remark} If an AZD exists on a line bundle $L$ on a smooth projective
variety $M$, $L$ is pseudoeffective by the condition 1 above.
\end{remark}
\begin{theorem}(\cite[Theorem 1.5]{tu,tu2,d-p-s})
Let $X$ be a smooth projective variety and let $L$ be a pseudoeffective 
line bundle on $X$.  Then $L$ has an AZD.
\end{theorem}
In fact an AZD of $L$ can be constructed as follows. 
Let $E$ be the set of singular hermitian metric on $L$ defined by
\[
E = \{ h ; h : \mbox{lowersemicontinuous singular hermitian metric on $L$}, 
\]
\[
\hspace{70mm}\Theta_{h}\,
\mbox{is positive}, \frac{h}{h_{0}}\geq 1 \}.
\]
We set 
\[
h_{min} = h_{0}\cdot\inf_{h\in E}\frac{h}{h_{0}},
\]
where the infimum is taken pointwise. 
Then it is easy to check that  $h_{min}$ is an AZD of $L$. 

The AZD $h_{min}$ constructed as above is said to be a canonical AZD of 
$L$. 
Also this construction can be generalized to the case of singular hermitian 
line bundles whose curvature current is bounded from below by some 
negative multiple of a K\"{a}hler form. 

\section{Inductive construction of metrics}

Let $X$ be a smooth projective variety and let 
$K_{X}$ be the canonical line bundle of $X$. 
Let $h$ be a canonical AZD of $K_{X}+L$.  
Let $n$ denote the dimension of $X$.

Let $A$ be a sufficiently ample line bundle on $X$ 
such that for every pseudoeffective singular hermitian 
line bundle $(F,h_{F})$
\[
{\cal O}_{X}(A+F)\otimes{\cal I}(h_{F})
\]
and 
\[
{\cal O}_{X}(K_{X}+L+A+F)\otimes{\cal I}(h_{F})\\
\]
are globally generated. 
This is possible by \cite[p. 667, Proposition 1]{si}. 
Let $h_{A}$ be a $C^{\infty}$ hermitian metric on $A$
 with strictly positive curvature. 

For $m\geq 0$, let $h_{m}$ be the singular hermitian metrics
on $A + m(K_{X}+L)$ constructed as follows. 
Let $h_{0}$ be a $C^{\infty}$-hermitian metric $h_{A}$ on 
$A$ with strictly positive curvature. 
Let $\{ \sigma_{0}^{(m)},\ldots ,\sigma_{N(m)}^{(m)}\}$
be an orthonormal basis of 
$H^{0}(X,{\cal O}_{X}(A+m(K_{X}+L))\otimes{\cal I}(h^{m-1}))$ 
with respect to the inner product :
\begin{eqnarray*}
(\sigma ,\sigma^{\prime}) & := &
\int_{X}\sigma\cdot \overline{\sigma^{\prime}}\cdot h_{A}\cdot h^{m}\cdot h_{L} \\
 &= & \int_{X}\sigma\cdot \overline{\sigma^{\prime}}\cdot (h_{A}\cdot h^{m}\otimes (dV)^{-1}\cdot h_{L})\cdot dV,   
\end{eqnarray*}
where $dV$ is an arbitrary nowhere degenerate $C^{\infty}$ volume form on $X$.
We set 
\[
K_{m} := \sum_{i=0}^{N(m)}\mid \sigma_{i}^{(m)}\mid^{2},
\]
where $\mid \sigma_{i}^{(m)}\mid^{2}$ 
denotes $\sigma_{i}^{(m)}\cdot \overline{\sigma_{i}^{(m)}}$.
We call $K_{m}$ the {\bf Bergman kernel} of $A + m(K_{X}+L)$ with respect 
to $h_{A}\cdot h^{m-1}\otimes (dV)^{-1}\cdot h_{L}$ and $dV$.
Clearly it is independent of the choice of the orthonormal basis.
And we define the singular hermitian metric $h_{m}$ 
on $A + m(K_{X}+L)$ by  
\[
h_{m} :=  K_{m}^{-1}.
\]
It is clear that $K_{m}$ has semipositive curvature 
in the sense of currents. 
We note that for every $x\in X$
\[
K_{m}(x) = \sup \{ \mid\sigma\mid^{2}(x) ;  
\sigma\in \Gamma (X,{\cal O}_{X}(A+mK_{X})), 
\int_{X}h_{A}h_{m}\cdot \mid\sigma\mid^{2} = 1\}
\]
holds by definition (cf. \cite[p.46, Proposition 1.4.16]{kr}).

Let $dV$ be a $C^{\infty}$-volume form on $X$ with respect to a 
K\"{a}hler form $\omega$ on $X$. 
For a singular hemitian line bundle $(F,h_{F})$ on $X$,
let $A^{2}(M,F,h_{F},dV)$ denote the Hilbert space of 
$L^{2}$ holomorphic sections of $F$ with respect to $h_{F}$ 
and $dV$. 
We may also assume that for any pseudoeffective
singular hemitian line bundle
$(F,h_{F})$ and for any point $x\in X$, 
there exists an interpolation operator
\[
I_{x} : A^{2}(x,K_{X}\otimes A\otimes F,dV^{-1}h_{A}h_{F},\delta_{x})
\longrightarrow  A^{2}(X,K_{X}\otimes A\otimes F,dV^{-1}h_{A}h_{F},dV)
\]
such that the operator norm of $I_{x}$ is bounded from above 
by a positive constant independent of  $x\in X$ and $(L,h_{L})$, 
where $\delta_{x}$ denotes the Dirac measure at $x$. 
This is certainly possible, if we take $A$ to be 
sufficiently ample.

In fact let $x$ be a point on $X$ and let $(U,z_{1},\ldots ,z_{n})$ be a local 
coordinate neighbourhood of $x$ which is biholomorphic to $\Delta^{n}$ 
and $z_{i}(x) = 0 (1\leq i\leq n)$. 
Then by  Theorem 2.1, we find an interpolation operator
\[
I_{x}^{U}: A^{2}(x,K_{X}\otimes A\otimes F,dV^{-1}h_{A}\cdot h_{F},\delta_{x})
\longrightarrow  A^{2}(U,K_{X}\otimes A\otimes F,dV^{-1}h_{A}\cdot h_{F},dV)
\]
such that the operator norm of $I_{x}^{U}$ is bounded from above 
by a positive constant $C_{U}$ independent of $(F,h_{F})$.
Now we note that the curvature of $h_{A}\cdot h_{F}$ is bounded from below by 
 the K\"{a}hler form $\Theta_{A}$.
Let $\rho$ be a $C^{\infty}$-function on $X$ such that 
$\mbox{Supp}\,\,\rho \subset\subset U$, $0\leq \rho\leq 1$ and 
$\rho \equiv 1$ on a neighbourhood of $x$. 
Let $\sigma_{x}$ be an element of $A^{2}(x,K_{X}\otimes A\otimes F,dV^{-1}h_{A}h_{L},\delta_{x})$.
Then  replacing  $(A,h_{A})$ by its  sufficiently high positive multiple,
we may assume that 
\[
\Theta_{A} + (n+1)\sqrt{-1}\partial\bar{\partial}(\rho\cdot\log \sum_{i=1}^{n}\mid z_{i}\mid^{2}) \geq \omega
\]
holds on $X$.  
We also note that there exists a positive constant $C^{\prime}_{U}$ 
independent of $(F,h_{F})$ and $\sigma_{x}$ such that 
\[
\int_{X}\exp (-(n+1)\rho\cdot\log \sum_{i=1}^{n}\mid z_{i}\mid^{2})\cdot\mid\bar{\partial}(\rho\cdot I^{U}_{x}\sigma_{x} )\mid^{2} dV \leq 
C^{\prime}_{U}\cdot (dV^{-1}\cdot h_{A}\cdot h_{F})(\sigma_{x},\sigma_{x})
\]
holds, where $\mid\bar{\partial}(\rho\cdot I^{U}_{x}\sigma_{x} )\mid^{2}$
denotes the norm with respect to $h_{A}\cdot h_{F}$ and $\omega$.
In fact $C^{\prime}_{U}$ only depends on $C_{U}$ and the supremum 
of the norm of $\bar{\partial}\rho$ with respect to $\omega$.
Then by the usual $L^{2}$-estimate, we may assume that we can solve the 
$\bar{\partial}$-equation 
\[
\bar{\partial}u = \bar{\partial}(\rho\cdot I^{U}_{x}\sigma_{x} )
\]
with 
\[
u (x) = 0
\]
so that 
\[
\int_{X}h_{A}\cdot h_{F}\mid u\mid^{2}
\leq C\cdot (dV^{-1}h_{A}\cdot h_{F})(\sigma_{x},\sigma_{x})
\]
holds for a positive constant $C$ independent of $(F,h_{F})$ and 
$\sigma_{x}$.
Then 
\[
\rho\cdot I^{U}_{x}\sigma_{x} - u
\in H^{0}(X,{\cal O}_{X}(K_{X}+A+L)\otimes{\cal I}(h_{F}))
\]
is an extension of $\sigma_{x}$. 
Since $X$ is compact, moving $x$ and $U$, by the above estimates this implies the assertion.
\begin{lemma}
Let $h$ be a canonical AZD of $K_{X}$ constructed as in the proof 
of Theorem 2.2.
Then the inclusion :
\[
{\cal I}(h^{m}) \subseteq {\cal I}(h_{m})
\]
holds for every $m\geq 0$.  
\end{lemma}

By the choice of $A$ and Lemma 3.1, $h_{m}$ is well defined 
for every $m\geq 0$.

Now we shall make the above lemma quantitative. 
\begin{lemma} 
There exists a positive constant  $C$ such that 
\[
h_{m} \leq C\cdot h_{A}\cdot h^{m} 
\]
holds for every $m\geq 0$. 
\end{lemma}
{\bf Proof}. 
Let us denote the Bergman kernel of 
$A +m(K_{X}+L)$ with respect to a singular hermitian metric 
$H$ on $A+m(K_{X}+L)$ and the volume form $dV$ by 
$K(A+m(K_{X}+L),H,dV)$. 
In this notation $K_{m}$ is expressed as 
$K(A+ m(K_{X}+L), h_{A}\cdot h^{m-1}\cdot h_{L}\otimes (dV)^{-1},dV)$.

Then by the $L^{2}$-extension form a point to $X$ as above we have that 
\[
K(A+ m(K_{X}+L),h_{A}\cdot h^{m-1}\cdot dV^{-1}\cdot   h_{L},dV) 
\geq  C^{-1}(h_{A}\cdot h^{m})^{-1}
\] 
hold for some positive constant $C$.  
This completes the proof of Lemma 3.2. 
{\bf Q.E.D.}
\vspace{5mm} \\
On the other hand by the submeanvalue property of plurisubhamonic functions, 
we have : 
\begin{lemma}
Let $x$ be a point on $M$ and let $\zeta = (\zeta_{1},\ldots , \zeta_{n})$ 
be a coordinate with center $x$. 
Then for every sufficiently small $r > 0$ 
\[
h_{A}\cdot K_{m}\leq \sup_{\mid\zeta\mid < r} h^{-m}\cdot \frac{1}{\pi^{n}r^{2n}/n!}
\]
hold, where we have trivialized $A$ and $K_{M}+L$ around $x$ by taking 
holomorphic frames. 
\end{lemma}
By Lemma 3.2 and Lemma 3.3 
\[
K_{\infty} : = \mbox{the uppersemicontinuous envelope of} 
\,\,\,\, \limsup_{m\rightarrow\infty}
\sqrt[m]{K_{m}}
\]
is a well defined volume form on $X$ which does not vanish 
outside of a set of measure $0$. 
We set 
\[
h_{\infty} : = \frac{1}{K_{\infty}}. 
\]
Then by Lemma 3.2, we see that 
\[
h_{\infty} \leq  h
\]
holds. 
By Lemma 3.3 we have the opposite estimate :  
\[
h_{\infty} \geq  h
\]
holds. 

Hence we have the following theorem.
\begin{theorem}
$h_{\infty} = h$  holds. 
In particular $h_{\infty}$ is an AZD of $K_{X} + L$.  
\end{theorem}
\section{Proof of Theorem 1.1}

Let $M, S, L$ be as in Theorem 1.1. Let $h_{S}$ be a canonical AZD of 
$K_{M}+L\mid_{S}$. 
Let $A$ be a sufficiently ample line bundle on $M$.  
Let us define the singular hemitian metric on $m(K_{M}+L)\mid_{S}$ by 
\[
h_{m,S} := K(A+m(K_{M}+L)\mid_{S},h_{A}\cdot h_{S}^{m-1}\cdot dV_{M}^{-1}\cdot 
h_{L},d\Psi_{S})^{-1}
\]
Then as in the last section, we see that 
\[
h_{\infty ,S} : = \liminf_{m\rightarrow\infty}\sqrt[m]{h_{m,S}}
\]
exists and an AZD of $K_{M}+L\mid_{S}$. 
We consider the Bergman kernel 
\[
K(S,A+m(K_{M}+L)\mid_{S},h_{A}\cdot h_{S}^{m-1}\cdot dV_{M}^{-1}\cdot 
h_{L},d\Psi_{S}) = \sum_{i}\mid \sigma_{i}^{(m)}\mid^{2}, 
\]
where $\{ \sigma_{i}^{(m)}\}$ is a complete orthonormal basis of 
$A^{2}(S,A+m(K_{M}+L)\mid_{S},h_{A}\cdot h_{S}^{m-1}\cdot dV_{M}^{-1}\cdot 
h_{L},d\Psi_{S})$. 
We note that (cf. \cite[p.46, Proposition 1.4.16]{kr})
\[
K(S,A+m(K_{M}+L)\mid_{S},h_{A}\cdot h_{S}^{m-1}\cdot dV_{M}^{-1}\cdot 
h_{L},d\Psi_{S})(x)
\]
\[= \sup \{ \mid\sigma\mid^{2}(x) \, \mid\, 
\sigma \in A^{2}(S,A+m(K_{M}+L)\mid_{S},h_{A}\cdot h_{S}^{m-1}\cdot dV_{M}^{-1}\cdot 
h_{L},d\Psi_{S}), \parallel \sigma\parallel = 1\}  
\]
holds for every $x\in S$.
We note that as in Lemma 3.2 there exists a positive constant $C_{0}$ independent of $m$ such that 
\[
h_{m,S} \leq C_{0}\cdot h_{A}\cdot h_{S}^{m}
\]
holds for every $m\geq 1$.

Inductively on $m$, we cextend each 
\[
\sigma\in
A^{2}(S,A+m(K_{M}+L)\mid_{S},h_{A}\cdot h_{S}^{m-1}\cdot dV_{M}^{-1}\cdot 
h_{L},d\Psi_{S}) 
\]
to a section 
\[
\tilde{\sigma}\in A^{2}(M,A+ m(K_{M}+L),dV^{-1}\cdot h_{L}\cdot \tilde{h}_{m-1},dV)
\]
with the estimate
\[
\parallel \tilde{\sigma}\parallel\leq C_{m} \parallel \sigma\parallel
\]
where  
$\parallel\,\,\,\,\parallel$'s denote the $L^{2}$-norms respectively, 
$C_{m}$ is the positive constant depending only on $m$ which will be specified later  
(the existence of $C_{m}$ is assured by Theorem 2.1), 
we have defined 
\[
\tilde{K}_{m}(x) := \sup\{ \mid \tilde{\sigma}\mid^{2}(x) \mid \,\,
\parallel\tilde{\sigma}\mid_{S}\parallel = 1, 
\parallel\tilde{\sigma}\parallel \leq C_{m} \}
\]
and set
\[
\tilde{h}_{m} = \frac{1}{\tilde{K}_{m}}.  
\]
Let us specify $C_{m}$. 
Let $X$ be the closed set as in Theorem 2.1. 
There exists a small Stein neighbourhood $W$ of $S\backslash X$ in $M$ 
and a holomorphic retraction 
\[
\rho : W \longrightarrow S\backslash X.
\]
Then $\rho$ defines a linear map 
\[
I_{\rho}^{m} : A^{2}(S,m(K_{M}+L)+A,h_{S}^{m-1}\otimes dV^{-1}\otimes h_{L}\otimes h_{A}) 
\]
\[
\hspace{30mm} \rightarrow   A^{2}(W,m(K_{M}+L)+A,\tilde{h}_{m-1}\otimes dV^{-1}\otimes h_{L}) 
\]
Then we may take $C_{m}$ as  
\[
C_{m}:= C\cdot \parallel I_{\rho}^{m}\parallel 
\]
where $C$ is a positive constant independent of $m$ 
and $\parallel I_{\rho}^{m}\parallel$ denotes the operator norm of $I_{\rho}^{m}$. 
This estimate follows from the proof of Theorem 2.1 (cf. \cite[p.9]{o}). 
By the above inductive estimates, we see that 
\[
\tilde{h}_{\infty} :=  \liminf_{m\rightarrow\infty}\sqrt[m]{\tilde{h}_{m}}
\]
exists and gives an extension of $h_{\infty ,S}$.
Then by Theorem 2.1 we may extend every element of 
$A^{2}(m(K_{M}+L)\mid_{S},dV_{M}^{-1}\cdot h_{L}\cdot \tilde{h}_{\infty}\mid_{S},dV_{M}[\Psi ])$ to 
$A^{2}(m(K_{M}+L),dV_{M}^{-1}\cdot h_{L}\cdot \tilde{h}_{\infty},dV_{M})$. 
This completes the proof of Theorem 1.1. 

\section{Generalization of Theorem 1.1}

Let $M$ be a smooth projective variety 
and let $(L,h_{L})$ be a singlar hermitian line bundle on $M$ such that 
$\Theta_{h_{L}}\geq 0$ on $M$.  
Let $dV$ be a $C^{\infty}$-volume form on $M$. 
Let $\sigma \in \Gamma (\bar{M},{\cal O}_{\bar{M}}(m_{0}L)\otimes {\cal I}(h))$ be a 
global section. 
Let $\alpha$ be a positive rational number $\leq 1$ and let $S$ be 
an irreducible subvariety of $M$ 
such that  $(M, \alpha (\sigma ))$ is logcanonical but not KLT(Kawamata log-terminal)
on the generic point of $S$ and $(M,(\alpha -\epsilon )(\sigma ))$ is KLT on the generic point of $S$ 
for every $0 < \epsilon << 1$. 
We set 
\[
\Psi = \alpha \log h_{L}(\sigma ,\sigma ).
\]
 
We shall assume that $S$ is not contained in the 
singular locus of $h$, where the singular locus of $h$ means the 
set of points where $h$ is $+\infty$. 

For the moment we shall consider the case that $S$ is smooth 
(when $S$ is not smooth, we just need to repeat the following 
argument after taking an embedded resolution of $S$). 
In this case $\Psi$ may not belong to $\sharp (S)$, 
since $\Psi$ may not have the prescribed singularity along $S$ 
as in the definition of $\sharp (S)$. 
But the proof of Theorem 2.1 and hence proof of Theorem 1.1
also works in this case except a minor difference. 
The difference is that 
$dV_{M}[\Psi ]$ (which is defined similarly as above) may have singularities along some Zariski closed subset of $S$. 
Let $d\mu_{S}$ be a $C^{\infty}$-volume form on $S$ and 
let $\varphi$ be the function on $S$ defined by
\[
\varphi := \log \frac{dV[\Psi ]}{d\mu_{S}}. 
\]
According to the singularity of $\varphi$, the proof of 
Theorem 1.1 must be modified as follows. 

Let $d$ be a positive integer such that 
$d > \alpha m_{0}$. 
Let $h_{S}$ be a canonical AZD of 
$(K_{M}+dL,e^{-\varphi}\cdot dV^{-1}\otimes h_{L}^{d})$, i.e., 
\[
h_{S} = \inf \{ H_{S} \mid 
\frac{H_{S}}{e^{-\varphi}\cdot dV^{-1}\otimes h_{L}^{d}}\geq 1, 
\Theta_{H_{S}}\geq 0\} ,
\]
where $H_{S}$ runs singular hermitian metrics on 
$K_{M}+ dL$ satisfying the above conditions. 
Let $A$ be a sufficiently ample line bundle on $M$. 
We set 
\[
K_{m} := K(A+m(K_{M}+dL)\mid_{S},h_{A}\cdot h_{S}^{m-1}
\cdot dV^{-1}\cdot h_{L}^{d},dV[\Psi ]). 
\]
Then by Theorem 2.1 (with the modified $dV[\Psi ]$), 
exactly as in the proof of Theorem 1.1,  inductively 
we may extend  
\[
h_{m}:= K_{m}^{-1}
\]
to a singular hermitian metric $\tilde{h}_{m}$ 
on $A + m(K_{M}+dL)$.  
Here we have used the factor $e^{-\varphi}$ to extend sections on $S$ 
inductively. 
Also by the same argument as in the proof of Theorem 1.1, 
\[
\tilde{h}_{\infty} := \liminf_{m\rightarrow\infty}\sqrt[m]{\tilde{h}_{m}}
\]
exists and is a singular hermitian metric on 
$K_{M}+dL$ with semipositive curvature.
Also as in Section 1, we see that 
$\tilde{h}_{\infty}\mid_{S}$ is an AZD of 
$(K_{M}+dL\mid_{S},e^{-\varphi}\cdot dV^{-1}\otimes h_{L}^{d})$.
Hence applying Theorem 2.1, we obtain the following theorem : 
\begin{theorem}
Let $M$,$S$,$\Psi$ be as above. 
Suppose that $S$ is smooth.  
Then every element of 
$A^{2}(S,{\cal O}_{S}(m(K_{M}+dL)),e^{-(m-1)\varphi}\cdot dV^{-m}\otimes 
h_{L}^{m},dV[\Psi ])$ 
extends to an element of 
\[
H^{0}(M,{\cal O}_{M}(m(K_{M}+dL))). 
\]
\end{theorem}
As we mentioned as above the smoothness assumption on $S$ is 
just to make the statement simpler. 

As an example of an application, we have : 
\begin{corollary}(\cite{t})
Let $\pi : X\longrightarrow \Delta$ be a semistable 
degeneration of projective variety over the unit disk. 
Let $X_{0} = \pi^{-1}(0) = \sum_{i}D_{i}$ 
be the irreducible decomposition. 
Then we have that 
\[
\sum_{i}P_{m}(D_{i}) \leq P_{m}(X_{t})
\]
holds where $t$ is any regular value  of $\pi$ 
and $P_{m}$ denotes the $m$-th plurigenus. 
\end{corollary}

Author's address\\
Hajime Tsuji\\
Department of Mathematics\\
Tokyo Institute of Technology\\
2-12-1 Ohokayama, Megro 152-8551\\
Japan \\
e-mail address: tsuji@math.titech.ac.jp

\begin{thebibliography}{99}
\bibitem{d-p-s}  J.P. Demailly-T. Peternell-M. Schneider : 
Pseudo-effective line bundles on compact K\"{a}hler manifolds, 
math. AG/0006025 (2000).
\bibitem{ka} Y. Kawamata, Subadjunction of log canonical divisors. II. Amer. J. Math. 120 (1998), no. 5, 893--899.  
\bibitem{kr} S. Krantz, Function theory of several complex variables, 
John Wiley and Sons (1982).  
\bibitem{o} T. Ohsawa, On the extension of $L^{2}$ holomorphic functions V,
effects of generalization, Nagoya Math. J. (2001) 1-21. 
\bibitem{o-t}T. Ohsawa and K. Takegoshi, $L^{2}$-extension of holomorphic
functions, Math. Z. 195 (1987),197-204.
\bibitem{si} Y.-T. Siu, Invariance of plurigenera, Invent. Math. 134 (1998), 661-673.
\bibitem{t} H. Tsuji, Deformation invariance of plurigenera, to appear. 
\bibitem{tu}H. Tsuji, Analytic Zariski decomposition, Proc. of Japan Acad.
61(1992) 161-163.
\bibitem{tu2} H. Tsuji, Existence and Applications of Analytic Zariski Decompositions, Analysis and Geometry in Several Complex Variables (Komatsu and Kuranishi ed.),  Trends in Math. 253-271, Birkh\"{a}user (1999).
\end{thebibliography}
\end{document}